\documentclass[12pt]{amsart}
\usepackage{amsfonts}
\usepackage{amssymb}
\newcommand{\dis}{\displaystyle}

\newcommand{\goes}{\rightarrow}
\newcommand{\Rr}{\mathbb{R}}
\newcommand{\Cc}{\mathbb{C}}

\newcommand{\Zz}{\mathbb{Z}}

\theoremstyle{plain}
\newtheorem{theorem}{Theorem}
\newtheorem{lemma}[theorem]{Lemma}

\newtheorem{corollary}[theorem]{Corollary}

\newcommand{\beq}{\begin{eqnarray}}
\newcommand{\eeq}{\end{eqnarray}}
\newcommand{\beqn}{\begin{eqnarray*}}
\newcommand{\eeqn}{\end{eqnarray*}}
\newcommand{\tr}[1]{\mbox{$\, ^t\! #1$}}
\newcommand{\Adj}{{\rm Adj}\,}

\begin{document}
\title{Which Singular K3 Surfaces Cover an Enriques Surface}
\author{Al\.i S\.inan Sert\"{o}z}
\address{B\.ilkent University \\ Department of Mathematics \\
TR-06533 Ankara, Turkey}
\email{sertoz@fen.bilkent.edu.tr}
\thanks{Research partially supported bt T\"{U}B{\.I}TAK-BDP}
\renewcommand{\subjclassname}{\textup{2000} Mathematics Subject Classification}
\date{May 2002}
\begin{abstract} 
We determine the necessary and sufficient conditions on the entries
of the intersection matrix of the transcendental lattice of a K3
surface for the K3 surface to doubly cover an Enriques surface.
\end{abstract}
\maketitle


\section{Introduction}

When $X$ is a singular K3 surface over the field $\Cc$, the transcendental
lattice $T_X$ of $X$ is denoted by its intersection matrix 
\beq
\left( \begin{matrix} 2a & c \\ c & 2b \end{matrix} \right) \label{eq:1}
\eeq
with respect to some basis $\{ u,v\}$, where $a,b >0$ and
$4ab-c^2>0$. For the definitions and basic facts about K3 surfaces we
refer to \cite{barth}.

If $U$ denotes the hyperbolic lattice of rank 2 and if $E_8$ denotes the
even unimodular negative definite lattice of rank 8, we then define a
sublattice $\Lambda^-$ of the K3-lattice $\Lambda$ as
\[ \Lambda^-=U\oplus U(2)\oplus E_8(2). \]
Following the works of Horikawa on the period map of Enriques surfaces
and work of Nikulin on the embeddings of even lattices, Keum gave an
integral lattice theoretical criterion for the existence of a fixed
point free involution on a K3 surface,
\cite{horikawa1,horikawa2,nikulin2, keum}. This criterion is then applied in
\cite{keum} to show that every Kummer surface is the double cover of
some Enriques surface, which corresponds to taking $a$, $b$, $c$ of
$T_X$ even and $17 \leq \rho (X)\leq 20$, see also \cite{nikulin1,morrison}. 

A K3 surface with $12\leq \rho (X) \leq 20$  covers
an Enriques surface if and only if there is a primitive embedding 
$\phi :T_X\goes\Lambda^-$ such that the orthogonal complement of the
image in $\Lambda^-$ contains no self intersection -2 vector, and when
$\rho(X)=10$ or $11$, one also needs to have $length\, (T_X)\leq
\rho(X)-2$, 
\cite[Theorem 1]{keum}.

We implement this criterion to find explicit
necessary and sufficient conditions on the entries of $T_X$ so that
$X$ covers an Enriques surface when $\rho (X)=20$. In practice, if $X$
actually covers an Enriques surface it is sometimes, but by no means
always, easy to exhibit an
embedding $\phi :T_X\goes\Lambda^-$ 
such that i) it is possible to
demonstrate that $\phi$ is primitive and that ii) it is possible to show
that the existence of a self intersection -2 vector in $\phi (X)^\perp$ 
leads to a contradiction. Moreover in case $X$ does not cover an
Enriques surface 
then it is hard work to demonstrate that for every embedding the
orthogonal complement of the image has a self intersection -2
vector. We resolve this difficulty in 
\begin{theorem} \label{thm:1}
If $X$ is a singular K3 surface with transcendental lattice given as in
(\ref{eq:1}), then $X$ covers an Enriques surface if and only if one of
the following conditions hold:\\
{\bf I }   $a$, $b$ and $c$ are even. (Keum's result, see \cite{keum}). \\
{\bf II}  $c$ is odd and $ab$ is even. \\
{\bf III-1 } $c$ is even. $a$ or $b$ is odd. The form $ax^2+cxy+by^2$
does not represent~$1$. \\
{\bf III-2 } $c$ is even. $a$ or $b$ is odd. The form $ax^2+cxy+by^2$
represents~$1$, and $4ab-c^2\neq 4,8,16$. 

Equivalently, $X$ fails to doubly cover an Enriques surface if and
only if one of the 
following conditions hold: \\
{\bf III-3 } $c$ is even. $a$ or $b$ is odd. The form $ax^2+cxy+by^2$
represents~$1$, and $4ab-c^2= 4,8,16$. \\ 
{\bf IV } $abc$ is odd. 
\end{theorem}


\section{Parities in Transcendental Lattice}
Before we proceed with the proof we must check that the parity
properties given in Theorem~\ref{thm:1} are well defined. 

Let $\dis \gamma=\left( \begin{matrix} x & y \\ z & w \end{matrix}
\right) \in
SL_2(\Zz )$. Then every matrix of the form $\tr{\gamma}T_X\gamma$
represents the transcendental lattice of $X$ with respect to some
basis. Setting 
\beqn
\tr{\gamma}\, T_X \, \gamma &=& \left( \begin{matrix}
2(ax^2+cxz+bz^2) & 2axy +c(xw+yz)+2bwz \\
2axy +c(xw+yz)+2bwz & 2(ay^2+cyw+bw^2)
\end{matrix} \right) \\
&=& \left( \begin{matrix} 2a' & c' \\ c' & 2b' \end{matrix} \right) ,
\eeqn
we see by inspection that \\
{\bf I } If $a$, $b$, and $c$ are even, then $a'$, $b'$ and $c'$ are
even. \\
{\bf II } If $c$ is odd and $ab$ is even, then $c'$ is odd and $a'b'$
is even. \\
{\bf III } If $c$ is even with $a$ or $b$  odd, then $c'$ is even with
$a'$ or $b'$  odd. \\
{\bf IV } If $abc$ is odd, then $a'b'c'$ is odd.

 
\section{Two Lemmas on Integral Lattices}
We require two lemmas on integral lattices. 
The first one is a
divisibility result which will enable us later to conclude that a
certain lattice is unimodular.
The second one is a
numerical implementation for the primitiveness of an embedding in
terms of the entries of the matrix of embedding. 

We freely use some fundamental concepts related to integral lattices
for which we can refer to \cite{barth,birkhof,degtyarev,milnor}.

Let $M=(\Zz^n,A)$ be an integral lattice where $A=\tr{A}$ is the
intersection matrix, and let $\alpha=(\alpha_1,\dots,\alpha_n)$ be a
primitive element, i.e. $\gcd (\alpha_1,\dots,\alpha_n)=1$. 
We denote by $<\alpha,\beta>_M$ the inner product of the vectors
$\alpha$ and $\beta$ in $M$.
Denote the
orthogonal complement of $\alpha$ in $M$ by $\alpha^\perp$. Let $\{
\beta_2,\dots,\beta_n\}$ be a basis of $\alpha^\perp$.

Since $\alpha^\perp$ is a primitive sublattice, its basis can be
extended to a basis of $M$, say by the addition of a vector $\beta_1$.
If $\beta_i=(b_{i1},\dots,b_{in})$, $i=1,\dots,n$, and 
$B=\left(
b_{ij}\right)$ is the $n\times n$ integral matrix formed by
the entries of the $\beta_i$'s, then $\det B=\pm 1$. 

Let $\Adj B =\left( m_{ij} \right)$. Then in particular we can write
\beq
\pm 1=\det B =b_{11}m_{11}+b_{12}m_{21}+\cdots+b_{1n}m_{n1},
\label{eq:det}
\eeq
where $m_{ij}$ is $(-1)^{i+j}$ times
the determinant of the $(n-1)\times (n-1)$ matrix
obtained from $B$ by deleting the $j^{th}$ row and $i^{th}$ column.  This
equation implies
\beq
\gcd (m_{11},m_{12},\dots,m_{n1})=1. \label{eqn:2}
\eeq

\begin{lemma} \label{lemma:1}
The index of $\alpha\oplus\alpha^\perp$ in $M$ divides
$<\alpha,\alpha>_M $.
\end{lemma}
\proof
Let $C$ be the $n\times n$ matrix obtained by
replacing the first row of $B$ by $\alpha$. Then $|\det C|$ is the index
of $\alpha\oplus\alpha^\perp$ in $M$. Let $D$ be the $(n-1)\times
(n-1)$ matrix defined by 
\beqn
CA\tr{C} &=& \left( 
\begin{array}{cccc}
<\alpha,\alpha>_M & 0 &  \dots & 0 \\
              0   &   &        &   \\
        \vdots    &   & D      &   \\
             0    &   &        &
\end{array} 
\right) .
\eeqn
Let
\beqn
X&=&\left( A\tr{C} \right)^1=\left( \begin{matrix} x_1 \\ \vdots \\
x_n \end{matrix} \right) , \\
Y&=& \left( \Adj C  \right)^1,
\eeqn
where the supscript denotes as usual the column number.
Then
\beqn
<\alpha,\alpha>_M \, Y &=& \det C \, X, \; \mbox{or equivalently} \\
x_i \, \det C &=& <\alpha,\alpha>_M \, m_{i1}, \; i=1,\dots,n.
\eeqn
If $\det C$ does not divide $<\alpha,\alpha>_M$, then there is a prime
factor $p$ of $\det C$ such that $p\nmid <\alpha,\alpha>_M$ and 
$p|m_{i1}$, $i=1,\dots,n$. This contradicts equation~\ref{eqn:2}. 
\qed 

Next we derive some working tools to recognize the primitiveness of an
embedding through the embedding matrix. 

Let $L_1$ and $L_2$ be two lattices
with base elements $e_1,...,e_n$ and $f_1,...,f_m$ respectively
where $m\geq n$. Assume that we have an embedding of $L_1$ into
$L_2$ given by
\[ \phi (e_i)=a_{i1}f_1+\cdots +a_{im}f_m,\;\; i=1,...,n \]
where the $a_{ij}$'s are integers. 

Set
\[ A=\left( a_{ij}\right)_{\substack{ 1\leq i\leq n \\ 1\leq j \leq m }} \]
and for any choice of integers $1\leq t_1,\dots,t_n\leq m$, define
\beqn
A(t_1,...,t_n) &=& \bigl( a_{it_j} \bigr)_{1\leq i,j\leq n} \\
\Delta (t_1,...,t_n) &=& \det A(t_1,...,t_n) \\
C(t_1,...,t_n) &=& \Adj \; A(t_1,...,t_n) \\
               &=& \bigl( c_{ij}(t_1,...,t_n) \bigr)_{1\leq i,j\leq n}
\eeqn

Let $z\in L_2$ be an element such
that $Nz\in \phi (L_1)$, for some positive integer $N$, 
i.e. there exist integers $c_1,...,c_n$
such that
\beqn
Nz&=&c_1\phi (e_1)+\cdots +c_n\phi (e_n).
\eeqn
We may assume without loss of generality that
\[ (N,c_1,...,c_n)=1. \]

\begin{lemma} With the above notation $\dis N|\Delta(t_1,...,t_n)$ for
every choice of integers $1\leq t_1,\cdots ,t_n\leq m$.
\end{lemma}
\proof
We have
\beq
Nz&=&\sum_{i=1}^n c_i\phi (e_i) \nonumber \\
&=&\sum_{i=1}^{n}c_i\bigl( \sum_{j=1}^{m}a_{ij}f_j\bigr) \nonumber \\
&=& \sum_{j=1}^{m}\bigl( \sum_{i=1}^{n}c_ia_{ij}\bigr) f_j \nonumber
\eeq
and
\[ N|\sum_{i=1}^{n}c_ia_{ij} \; \mbox{for every $j=1,...,m$.} \]
In particular for any choice of integers
$d_1,...,d_m$ we have
\beqn
N&|& \sum_{j=1}^{m} d_j \bigl( \sum_{i=1}^{n}c_ia_{ij}\bigr) 
\eeqn
or equivalently
\beq
N &|& \sum_{i=1}^{n}c_i\bigl(
\sum_{j=1}^{m}d_ja_{ij}\bigr). \label{eq:3}
\eeq

For any choice of integers $1\leq t_1,\cdots ,t_n\leq m$ and $1\leq
s\leq n$ define 
for $j=1,\dots,m$
\beqn
d_j &=& 
\begin{cases}
c_{ks}(t_1,...,t_n) & {\rm if~}j=t_k, \\
0                   & {\rm otherwise}
\end{cases}
\eeqn
With these choices of $d_j$'s 
it follows that for $i=1,...,n$, we have
\beqn
\sum_{j=1}^{m}d_ja_{ij} &=& \sum_{k=1}^{n}a_{ik}c_{ks}(t_1,\dots,t_n) \\
&=& \delta_{is}\Delta(t_1,...,t_n),
\eeqn
and equation (\ref{eq:3}) becomes
\beqn
N&|& c_s \Delta(t_1,...,t_n) \; \mbox{for every $s=1,...,n$.} 
\eeqn
Since $(N,c_1,...,c_n)=1$, it follows that $N|\Delta(t_1,...,t_n)$ as required.
\qed 

Set
\beqn
d&=&gcd\{ \Delta(t_1,...,t_n) \; |\; 1\leq t_1,\cdots ,t_n\leq m\; \} .
\eeqn

\begin{corollary} If $d=1$, then $\phi$ is primitive.
\qed
\end{corollary}

\begin{lemma} If $d>1$, then the embedding $\phi$ is not primitive.
\end{lemma} 
\proof
When $d>1$, we explicitly construct  an element $z\in
L_2$ and an integer $N>1$ such that $z\notin \phi(L_1)$ but
$Nz\in\phi(L_1)$. 

Let $d=p_1^{\alpha_1}\cdots p_r^{\alpha_r}$ be the prime
decomposition of $d$. 
Without loss of generality assume, after
rearranging basis if necessary, that $\Delta(1,...,n)=dK$ where
$p_1\cdots p_r \nmid K$. For convenience set the notation as 
\beqn
A_0 &=& A(1,...,n)=\bigl( a_{ij} \bigr)_{1\leq i,j\leq n} \\
\Delta_0 &=& \Delta(1,...,n) \\
C_0 &=& C(1,...,n)=\bigl( c_{ij} \bigr)_{1\leq i,j\leq n} = \Adj A_0
\eeqn
Furthermore define
\beqn
d_i &=& gcd(c_{i1},...,c_{in}),\;\; i=1,...,n. 
\eeqn

Assume that $d|d_i$ for all $i$. Then $d^n|\Delta_0^{n-1}$, and 
$d|K^{n-1}$ which contradicts the choice of $\Delta_0$. 

Thus $d\nmid d_\ell$ for some $1\leq \ell \leq n$.

Set $d_0=gcd(d,d_\ell )$ and note that $1\leq d_0<d$.

Define the integers $k_1,...,k_m$, $c_1,...,c_n$ and $N$
as
\beqn
k_i &=& \left\{ \begin{array}{ll}
        0&{\rm if~} i\in \{1,...,\hat{\ell},...,n\} \\
        \Delta_0/d_0 & {\rm if~} i=\ell \\
        \Delta(1,\dots,\ell-1,i,\ell+1,\dots,n)/d_0 &{\rm if~} i=n+1,...,m.
        \end{array} \right. \\
c_i &=& c_{\ell i}/d_0, \;\; i=1,...,n \\
N &=& d/d_0.
\eeqn
We claim that
\beq
N &|& k_i,\;\; i=1,...,m, \label{eq:5} \\
gcd(N,c_1,...,c_n)&=&1, \label{eq:6} \\
{\rm and} && \nonumber \\
(c_1,...,c_n)A &=&(k_1,...,k_m). \label{eq:7}
\eeq
To prove (\ref{eq:5}), note that $d|\Delta(t_1,...,t_n)$ for any choice
of integers $1\leq t_1,...,t_n\leq m$, so
$N=(d/d_0)|(\Delta(t_1,...,t_n)/d_0)$. From the definition of the $k_i$'s it
follows now that $N|k_i$ for all $i=1,...,m$.

To prove (\ref{eq:6}) recall that $gcd(c_{\ell 1},...,c_{\ell
n})=d_{\ell}$ and  $d_0=gcd(d,d_{\ell})$ so
$gcd(N=d/d_0, c_1=c_{\ell 1}/d_0,...,c_n=c_{\ell n}/d_0)=1$.

To prove (\ref{eq:7}) it suffices to observe that
\beqn
(c_{\ell 1},...,c_{\ell n})\left( \begin{array}{c}a_{1i}\\
\vdots \\ a_{ni} \end{array} \right) &=& 
\begin{cases}
\delta_{\ell i}\Delta_0 & {\rm ~if}\;\; i=1,...,n \\
\Delta(1,\dots,\ell-1,i,\ell+1,\dots,n) & {\rm if~} i=n+1,...,m.
\end{cases}
\eeqn

Finally we put these together to define the torsion element $z$ as
\beqn
z&=& \frac{k_1}{N}f_1+\cdots+\frac{k_m}{N}f_m\in L_2.
\eeqn
It is now clear that 
\beqn
Nz&=& c_1\phi (e_1)+\cdots+c_n\phi (e_n)\in \phi (L_1),
\eeqn
and
\beqn
z&\not\in&\phi (L_1).
\eeqn
Hence we conclude that if
$d>1$, then $\phi$ is not primitive. \qed

The following theorem follows as a corollary of these two lemmas.

\begin{theorem}\label{thm:prm}
A lattice embedding is primitive if and only if the greatest common
divisor of the maximal minors of the embedding matrix
with respect to any choice of basis is $1$.
\qed
\end{theorem}

Though we gave a constructive proof of this theorem it is possible to
derive it from general principles. The rank of the matrix $A$ being
$n$, it has an integral diagonal form where the principal diagonal
entries are $h_1,\dots,h_n$, and all other entries are $0$. If $g_i$
denotes the greatest common divisor of all $i\times i$ minors of $A$,
then $g_i|g_{i+1}$ and $g_i=h_1\cdots h_i$.
In particular $g_n=1$ if and only if all the $h_i$'s
are $1$, i.e. $\phi:L_1\goes L_2$ is primitive if and only if there
exists a basis $\{e_1,\dots,e_n\}$ of $L_1$ and a basis
$\{f_1,\dots,f_m\}$ of $L_2$ such that $\phi(e_i)=f_i$,
$i=1,\dots,n$. It follows that $g_n$, the greatest common divisor of
the maximal minors of the embedding matrix $A$, is the order of
torsion of the quotient $L_2/\phi(L_1)$.

As an immediate application of this theorem we can indicate that all
the mappings in \cite[pp106-108]{keum} have embedding matrices whose
maximal minors have greatest common divisor equal to 1.


\section{The case when $c$ is even with $a$ or $b$ odd}

If $a$ is even, then set 
$\dis \gamma=\left( \begin{matrix} 2 & 1 \\ 1 & 1 \end{matrix} \right)$. 
If $\dis \tr{\gamma} \, T_X \, \gamma= 
\left( \begin{matrix} 2a' & c' \\
c' & 2b' \end{matrix} \right)$, then $a'$ and $b'$ are odd, and $c'$
is even. If $b$ is even, then  $\dis \gamma=\left( 
\begin{matrix} 1 & 1 \\ 1 & 2 \end{matrix} \right)$ changes $T_X$ into 
an  equivalent form where
again $a'$ and $b'$ are odd, and $c'$ even. So we might assume without
loss of generality that $ab$ is odd, and $c$ is even. 

We will consider a particular embedding of $T_X$ into
$\Lambda^-=U\oplus U(2)\oplus E_8(2)$.

Let $\{ u,v\}$ be a basis of $T_X$, 
$\{ u_1, u_2\}$ be a basis of  $U$ and $\{ v_1,v_2\}$ be a basis
of $U(2)$. 

Define $\phi:T_X\goes\Lambda^-$ by 
\beqn
\phi (u) &=& u_1+au_2, \\
\phi(v) &=& u_1+(c-a)u_2 +v_1+\frac{1}{2}(a+b-c)v_2.
\eeqn
It can be shown by direct computation that this is an embedding and by
theorem \ref{thm:prm} that this embedding is primitive.


\subsection{The form $ax^2+cxy+by^2$ does not represent $1$}
\mbox{} \\
Let $f=xu_1+x'u_2+yv_1+y'v_2+e\in \Lambda^-$, where $e\in E_8(2)$ with
$e\cdot e=-4k, \; k\geq 0$. (we will use $\cdot$ to denote the inner
product on $\Lambda^-$).

Impose the condition that $f$ lies in the orthogonal complement of
$\phi \left( T_X\right)$ in $\Lambda^-$ and that $f\cdot f=-2$.

Solving the equations $f\cdot \phi (u)=0$, $f\cdot \phi (v)=0$ for
$x'$, $y'$ and substituting into the equation $f\cdot f=-2$ gives
\beq
1-(ax^2+(c-2a)xy+(a+b-c)y^2)=2k \geq 0. \label{eq:8}
\eeq
The binary quadratic form $ax^2+(c-2a)xy+(a+b-c)y^2$ is equivalent to
the form $ax^2+cxy+by^2$. Since $a>0$ and $c^2-4ab<0$, this is a
positive definite form. Equation (\ref{eq:8}) holds if and only if
this form represents $1$, and then $k=0$. (see \cite{niven})

If we assume that the form $ax^2+cxy+by^2$ does not represent $1$,
then equation (\ref{eq:8}) cannot be solved, so there is no self
intersection $-2$ vector in the orthogonal complement of $\phi \left(
T_X\right)$.

This proves {\bf III-1}.


\subsection{The form $ax^2+cxy+by^2$ does represent $1$}
\mbox{} \\
In this case the binary quadratic form $ax^2+cxy+by^2$ is equaivalent
to the form $x^2+(ab-c^2/4)y^2$, see
\cite[p174]{niven}. Then a basis $\{ u, v\}$ of the transcendental
lattice exists such that with respect to that basis the matrix 
\beqn
T_X &=& \left( 
\begin{matrix} 2 (1) & 0 \\ 0 & 2 (\frac{\Delta}{4})  \end{matrix}
\right) 
\eeqn
where $\Delta=4ab-c^2$.

Let $\phi$ be a primitive embedding of $T_X$ into $\Lambda^-$ and 
set $\phi (u)=\alpha$ with
\[ \alpha=a_1u_1+a_2u_2+a_3v_1+a_4v_4+\omega_1 \]
where $\omega_1\in E_8(2)$ with $\omega\cdot\omega=-4k\leq 0$.

$\alpha\cdot\alpha=2$ forces $a_1$ and $a_2$ to be odd.

If $\beta=b_1u_1+b_2u_2+b_3v_1+b_4v_4+\omega_2$ is in the orthogonal
complement $\alpha^\perp$ of $\alpha$ in $\Lambda^-$, then
$\beta\cdot\alpha=0$ forces $b_1$ and $b_2$ to be of the same
parity. This in turn implies the following
\begin{lemma}\label{lemma:son}
If $\beta,\gamma\in \alpha^\perp$,  then $\beta\cdot\gamma \equiv 0
\mod 2$.
\end{lemma}

Let $\beta_1,\dots,\beta_{11}$ be basis elements for
$\alpha^\perp$, and $\dis B'=\left( 2b_{ij}\right)$,
$2b_{ij}=\beta_i\cdot\beta_j$ the intersection matrix for this
basis. Set $B=\left( b_{ij}\right)$. 

Let $C$ be the $12\times 12$-matrix whose rows are the coordinates of
$\alpha,\beta_1,\dots,\beta_{11}$ with respect to the standard basis
of $\Lambda^-$. And finally let $A$ denote the intersection matrix of
$\Lambda^-$ with respect to its standard basis. We have
\beq
CA\tr{C} &=& \left( 
\begin{array}{cccc}
2 & 0 &  \dots & 0 \\
              0   &   &        &   \\
        \vdots    &   & B'      &   \\
             0    &   &        &
\end{array} 
\right) . \label{eq:9}
\eeq
Since $\alpha,\beta_1,\dots,\beta_{11}$ is not a basis of $\Lambda^-$,
$|\det C|>1$. By lemma \ref{lemma:1}, $|\det C|$ divides $2$, hence is
equal to $2$. By interchanging $\beta_1$ by $\beta_2$ if necessary, we
can assume without loss of generality that $\det C=2$. 

It then follows from equation (\ref{eq:9}) that $\det B=1$.

Define a new lattice $L=\left( \Zz^{11}, B(-1)\right)$. $L$ has
signature $(\tau^+,\tau^-)=(10,1)$. Since $\tau^+-\tau^-\not\equiv 0
\mod 8$, $L$ is odd. Then $L$ is an indeterminate, odd, unimodular
lattice, and as such is isomorphic to $<-1>^1\oplus <1>^{10}$. 

There is an isomorphism $F:\alpha^\perp\goes L$ which sends $\beta_i$
to $e_i=(0,\dots,1,\dots,0)$, where $1$ is in the $i$-th place. This
isomorphis respects inner products in the sense that
\[-2[F(\lambda_1)\cdot F(\lambda_2)]=\lambda_1\cdot\lambda_2, \; 
\mbox{for all $\lambda_1,\lambda_2\in \alpha^\perp$}. \]

Let $e_1',\dots,e_{11}'$  be a basis of $L$ diagonalizing its
intersection matrix. Then the intersection matrix of
$\alpha\oplus\alpha^\perp$ with respect to the basis  $\alpha,
F^{-1}(e_1'),\dots,F^{-1}(e_{11}')$ is 
\beqn
\left( \begin{matrix} 
2 & 0   & 0  & \dots & 0 \\
0 & 2   & 0  & \dots & 0 \\
0 & 0   & -2 & \dots & 0 \\
\vdots  & \vdots & \vdots & \ddots & \vdots \\
0 & 0   & 0 & \dots & -2 
\end{matrix} \right)
\eeqn

We are looking for the existence of a primitive
embedding 
\[ \phi : T_X\longrightarrow \alpha\oplus\alpha^\perp\subset \Lambda^-
\]
such that with respect to this new basis of $\alpha\oplus\alpha^\perp$, 
\beqn
\phi(u) &=& (1,0,\dots,0), \\
\phi(v) &=& (0,x_0,\dots,x_{10})
\eeqn
such that
\beqn
\phi(v)\cdot\phi(v)=2x_0^2-2x_1^2-\cdots-2x_{10}^2=2 \left(
\frac{\Delta}{4} \right) .
\eeqn

Using theorem \ref{thm:prm}, the problem reduces to a problem in the
lattice $L$, that of investigating the existence of integers $x_0,\dots,
x_{10}$ such that if $x=(x_0,\dots,x_{10})\in L$ then the following
conditions are satisfied:
\beq
gcd(x_0,\dots,x_{10})&=&1, \nonumber \\
x\cdot x &=& -x_0^2+x_1^2+\cdots+x_{10}^2 \nonumber \\
&=& -\left( \frac{\Delta}{4} \right) , \; {\rm and} \nonumber \\
y\cdot x =0 &\Longrightarrow & y\cdot y \neq 1, \; \mbox{for every
$y\in L$}. \label{condition}
\eeq
The existence of such integers is equivalent to $X$ covering an
Enriques surface.

The set of negative self intersection elements of $L$ span an open
convex cone in $\Zz^{11}\otimes \Rr$, and we refer to \cite{benedetti} for
details. We will utilize the techniques of Vinberg from \cite{vinberg}
to understand the existence of integers as above.

All automorphisms of $L$ are generated by reflections and a
fundamental region for negative self intersecting vectors in $L$
is bounded by reflecting
hyperplanes. Studying the nature of these hyperplanes, as in
\cite{vinberg}, we conclude that the conditions in the set of equations
(\ref{condition}) holds if and only if 
\beqn
gcd(x_0,\dots,x_{10})=1, \\
-x_0^2+x_1^2+\cdots+x_{10}^2 = -\left( \frac{\Delta}{4} \right), \\
x_1\geq \cdots \geq x_{10}>0, \\
x_0 \geq x_1+x_2+x_3, \; {\rm and} \\
3x_0 > x_1+\cdots+x_{10}.
\eeqn

Let $P$ denote the set of all $x\in L$ satisfying the above
conditions.

The rest of this case is elementary and we summarize the results in
two technical lemmas.

\begin{lemma}
There is no $x\in P$ with $x\cdot x=-1,-2,-4$.
\end{lemma}
\proof Let $P(m)=\{x\in P\; | \; x=(m,x_1,\dots,x_{10}) \; \}$. 
Then it is easy to show that \\
$\dis \max_{x\in P(3m)}(x\cdot x)=5-4m$, $m>2$, \\
$\dis \max_{x\in P(6)}(x\cdot x)=-5$, \\
$\dis \max_{x\in P(3m+1)}(x\cdot x)=1-4m$, $m\geq 1$, \\
$\dis \max_{x\in P(3m+2)}(x\cdot x)=9-8m$, $m\geq 3$, \\
$\dis \max_{x\in P(8)}(x\cdot x)=-12$, \\
$\dis \max_{x\in P(5)}(x\cdot x)=-7$. \\ \mbox{} \\
These maximum values are achieved by the vectors \\
$[3m,m,\dots,m,m-2,1]$, \\
$[6,2,\dots,2,1,1,1]$, \\
$[3m+1,m+1,m,\dots,m,1]$, \\
$[3m+2,m+2,m,\dots,m,3]$, \\
$[8,4,2,\dots,2]$, \\
$[5,3,1,\dots,1]$, respectively.

It is now clear that $-1$ and $-2$ are never achieved. 
And if $x\cdot x=-4$, then $x\in P(4)$. But none of the vectors in
$P(4)$ achieve $-4$.
\qed

\begin{lemma}
For every positive integer $N$, other than $1,2$ and $4$, there is an
$x\in P$ such that $x\cdot x=-N$.
\end{lemma}
\proof 
In the following table each of the given vectors is in $P$,
and moreover 
$X_m(k)\cdot X_m(k)=-(m+24k)$, $Y_m\cdot Y_m=-m$, $Z(n)\cdot
Z(n)=-(4n-1)$, and $W(n)\cdot W(n)=-(4n-3)$. This then proves the
lemma.

$X_0(k)=[9k+4,3k+2,3k+1,\dots,3k+1,3k-1,2]$, $k\geq 1$. \\
$X_2(k)=[9k+4,3k+2,3k+1,\dots,3k+1,3k,3k,2]$, $k\geq 1$. \\
$X_4(k)=[12k+4,4k+2,4k+1,\dots,4k+1,4k,1]$, $k\geq 1$. \\
$X_6(k)=[9k+5,3k+2,3k+2,3k+1,\dots,3k+1,2]$, $k\geq 1$. \\
$X_8(k)=[9k+7,3k+3,3k+2,\dots,3k+2,3k,2]$, $k\geq 1$. \\
$X_{10}(k)=[12k+7, 4k+3, 4k+2, \dots, 4k+2, 4k+1,1]$, $k\geq 0$. \\
$X_{12}(k)=[12k+9, 4k+3, \dots, 4k+3,4k+2, 4k+1,1]$, $k\geq 0$. \\
$X_{14}(k)=[9k+8,3k+3,3k+3,3k+2,\dots,3k+2,2]$, $k\geq 0$. \\
$X_{16}(k)=[9k+10,3k+4,3k+3,\dots,3k+3,3k+1,2]$, $k\geq 1$. \\
$X_{18}(k)=[12k+12, 4k+4, \dots, 4k+4,4k+3, 4k+2,1]$, $k\geq 0$. \\
$X_{20}(k)=[6k+12, 2k+6,2k+3,\dots,2k+3,4]$, $k\geq 1$. \\
$X_{22}(k)=[9k+11,3k+4,3k+4,3k+3,\dots,3k+3,2]$, $k\geq 0$. \\
$Y_6=[4,1,\dots,1]$. \\
$Y_8=[6,2,\dots,2,1,1,1,1]$. \\
$Y_{20}=[6,2,2,1,\dots,1]$. \\
$Z(n)=[3n+1,n+1,n,\dots,n,1]$, $n\geq 1$. \\
$W(n)=[3n,n,\dots,n,n-1,n-1,1]$, $n\geq 2$.
\qed

Recalling that $\Delta=4ab-c^2$, these two lemmas complete the proofs
of {\bf III-2} and {\bf III-3}.


\section{The other cases}
Let $\{ u, v\}$ be a basis of the trancendental lattice giving the
matrix representation as in (\ref{eq:1}), and as before let $\{
u_1,u_2\}$ be the basis of $U$, and $\{ v_1,v_2\}$ the basis of $U(2)$.


\subsection{$abc$ is odd}
\mbox{} \\
Consider the mapping $\phi :T_X\goes \Lambda^-$ defined generically as
\beqn
\phi(u)&=& a_1u_1+a_2u_2+a_3v_1+a_4v_2+\omega_1 \\
\phi(v)&=& b_1u_1+b_2u_2+b_3v_1+b_4v_2+\omega_2 
\eeqn
where the $a_i$'s and $b_i$'s are integers, $\omega_i\in E_8(2)$. If
$\phi(u)\cdot\phi(u)=2a$ and 
$\phi(v)\cdot\phi(v)=2b$, then $a_1,a_2,b_1$ and $b_2$ are odd. But
this forces $\phi(u)\cdot\phi(v)$ to be even. Hence $T_X$ has no
embedding into $\Lambda^-$.

This proves {\bf IV}.  


\subsection{$c$ is odd and $ab$ is even}
\mbox{} \\
Consider the mapping $\phi :T_X\goes \Lambda^-$ defined as
\beqn
\phi(u) &=& au_1+u_2+\frac{1}{2}(c-ab-1)v_1, \\
\phi(v) &=& u_1+bu_2+v_2.
\eeqn
This is an embedding and by theorem \ref{thm:prm} it is primitive. Let
\[ f=c_1u_1+c_2u_2+c_3v_1+c_4v_2+\omega\in \Lambda^- \]
where $\omega\in E_8(2)$. 

$f\cdot\phi(u)=0$ and $f\cdot\phi(v)=0$ forces $c_1c_2$ to be
even. Then $f\cdot f\equiv 0 \mod 4$ and hence cannot be $-2$.

This proves {\bf II} and completes the proof of theorem \ref{thm:1}.


{\bf Acknowledgements: } I thank my colleagues A. Degtyarev,
A. Kerimov, A. Klyachko and E. Yal\c{c}{\i}n for numerous
discussions. 


\end{document}